# Orbit equivalence rigidity

By ALEX FURMAN

## Abstract

Consider a countable group $\Gamma$ acting ergodically by measure preserving transformations on a probability space $(X, \mu)$, and let $\mathcal{R}_\Gamma$ be the corresponding orbit equivalence relation on $X$. The following rigidity phenomenon is shown: there exist group actions such that the equivalence relation $\mathcal{R}_\Gamma$ on $X$ determines the group $\Gamma$ and the action $(X, \mu, \Gamma)$ uniquely, up to finite groups. The natural action of $\mathrm{SL}_n(\mathbb{Z})$ on the $n$-torus $\mathbb{R}^n/\mathbb{Z}^n$, for $n > 2$, is one of such examples. The interpretation of these results in the context of von Neumann algebras provides some support to the conjecture of Connes on rigidity of group algebras for groups with property T. Our rigidity results also give examples of countable equivalence relations of type $\mathrm{II}_1$, which *cannot* be generated (mod 0) by a *free action* of any group. This gives a negative answer to a long standing problem of Feldman and Moore.

## 1. Introduction and statements of main results

This paper contains applications of the rigidity result in [Fu] to some problems on measurable equivalence relations of type $\mathrm{II}_1$.

Let $(X, \mathcal{B}, \mu)$ be a nonatomic probability measure space, where $(X, \mathcal{B})$ is a standard Borel space. Given a measure preserving, ergodic action of a countable group $\Gamma$ on $(X, \mathcal{B}, \mu)$, one can define an equivalence relation $\mathcal{R}_\Gamma$ measurable with respect to $\mathcal{B} \times \mathcal{B}$ on $X \times X$ defined by

$$\mathcal{R}_\Gamma = \{(x, y) \in X \times X \mid \Gamma \cdot x = \Gamma \cdot y\}.$$

An especially interesting situation is that of an *essentially free* action $(X, \mu, \Gamma)$, i.e. an action which is free on an invariant subset $X' \subseteq X$ of full $\mu$-measure. We shall discuss the following basic

*Question* I. Let $(X, \mu, \Gamma)$ be an essentially free, finite measure preserving, ergodic action. Which properties of the group $\Gamma$ and the action $(X, \mu, \Gamma)$ are determined by the equivalence relation $\mathcal{R}_\Gamma$ on $(X, \mu)$?

More precisely, suppose $(X, \mu, \Gamma)$ and $(Y, \nu, \Lambda)$ are essentially free ergodic finite measure preserving (f.m.p. for short) actions which are *orbit equivalent*, i.e. that there exists a Borel isomorphism $\theta : X' \to Y'$ of full measure subsets, such that

$$\theta_* \mu \sim \nu, \quad \text{and} \quad (x_1, x_2) \in \mathcal{R}_\Gamma \iff (\theta\, x_1, \theta\, x_2) \in \mathcal{R}_\Lambda .$$

The question is what is the relationship between the groups $\Gamma$ and $\Lambda$ and the the actions $(X, \mu, \Gamma)$ and $(Y, \nu, \Lambda)$.

In this paper we establish an *orbit equivalence rigidity* of certain f.m.p. actions, namely we shall establish a collection of examples where the group $\Gamma$ and the action $(X, \mu, \Gamma)$ are completely determined (up to finite groups) just by the relation $\mathcal{R}_\Gamma$. Prior to the presentation of these results, let us recall some main known facts related to Question I.

Classical, by now, results of Dye [Dy], Ornstein and Weiss [OW], and for the most general case Connes, Feldman and Weiss [CFW], assert that all actions of (say countable) amenable groups are orbit equivalent to each other. On the other hand, an essentially free f.m.p. action of a nonamenable group is never orbit equivalent to an action of an amenable one (cf. Zimmer [Zi1], or [Zi4, 4.3.3]). Hence, amenability of an essentially freely acting group is determined by $\mathcal{R}_\Gamma$, and if the group is indeed amenable, no other information is contained in $\mathcal{R}_\Gamma$.

Another general property of $\Gamma$ which is determined by the orbit relation $\mathcal{R}_\Gamma$, generated by any essentially free f.m.p. $\Gamma$-action, is Kazhdan's property T (see [Fu, 1.4] and Zimmer [Zi4, 9.1.7(b)]).

Focusing on a specific class of groups, namely lattices in (semi)simple Lie groups, Zimmer obtained a sequence of rigidity results for actions of these groups. The main tool in these results is superrigidity for measurable cocycles (Zimmer [Zi2]), which is a generalization of Margulis's superrigidity for higher rank lattices (we refer to [Zi4] for the general exposition of these ideas). In particular, consider a lattice $\Gamma$ in a higher rank simple Lie group $G$ and a lattice $\Lambda$ in a (semi)simple Lie groups $H$, where $G$ and $H$ have finite centers. Suppose that $\Gamma$ and $\Lambda$ have orbit equivalent essentially free f.m.p. actions $(X, \Gamma)$ and $(Y, \Lambda)$. Then the ambient groups $G$ and $H$ are locally isomorphic, and the induced $G$-actions on $G/\Gamma \times X$ and $G/\Lambda \times Y$ are virtually isomorphic ([Zi2] or [Zi4] 4.5.2.1, see also Theorem E). Further, Zimmer showed that the assumption "$\Lambda$ is a lattice in a semisimple $H$," can be replaced by "$\Lambda$ is a countable group which admits a *linear representation* with an infinite image", the conclusion being that $\Lambda$ is virtually isomorphic to a lattice $\Lambda_1$ in $\operatorname{Ad} G$, and $(G/\Gamma \times X, G)$ is virtually isomorphic to $(\operatorname{Ad} G/\Lambda_1 \times Y_1, \operatorname{Ad} G)$ (see [Zi3]). Another fact, related to our discussion, is the result by Stuck and Zimmer [SZ], which states that any f.m.p. nonatomic ergodic action of a higher rank lattice is essentially free, up to a division by the finite center.

In the sequel we use a notion of *virtual isomorphism* which basically means isomorphism up to finite groups:

*Definition* 1.1. We shall say that two countable groups $\Gamma_1$ and $\Gamma_2$ are *virtually isomorphic groups* if there exist finite normal subgroups $N_i \triangleleft \Gamma_i$, $i = 1, 2$, such that the quotient groups $\Gamma'_i = \Gamma_i/N_i$ contain isomorphic subgroups of finite index: $\Gamma''_1 \cong \Gamma''_2$ where $[\Gamma'_i : \Gamma''_i] < \infty$. Let $(X_i, \mu_i, \Gamma_i)$, $i = 1, 2$ be two (ergodic) f.m.p. actions of virtually isomorphic groups $\Gamma_1, \Gamma_2$ as above. Then $\Gamma'_i = \Gamma_i/N_i$ naturally acts (ergodic f.m.p. action) on the quotient probability space $(X'_i, \mu'_i) = (X_i, \mu_i)/N_i$. Consider the action of the finite index subgroup $\Gamma''_i \subseteq \Gamma'_i$ on $(X'_i, \mu'_i)$ and let $(X''_i, \mu''_i)$ be one of at most $[\Gamma'_i : \Gamma''_i]$-many mutually isomorphic $\Gamma''_i$-ergodic components. We shall say that $(X_1, \mu_1, \Gamma_1)$ and $(X_2, \mu_2, \Gamma_2)$ are *virtually isomorphic actions* if $\Gamma_1, \Gamma_2$ are virtually isomorphic groups, and for some choice of $N_i$, $\Gamma''_i$ as above, the actions $(X''_1, \mu''_1, \Gamma''_1)$ and $(X''_2, \mu''_2, \Gamma''_2)$ are isomorphic.

We shall also need the following weaker version of orbit equivalence (see also Definition 2.2): two f.m.p. actions $(X, \mu, \Gamma)$ and $(Y, \nu, \Lambda)$ are *weakly orbit equivalent* if there exist positive measure subsets $A \subseteq X$ and $B \subseteq Y$, such that the restricted relations $\mathcal{R}_\Gamma|_A := \mathcal{R}_\Gamma \cap A \times A$ and $\mathcal{R}_\Lambda|_B := \mathcal{R}_\Lambda \cap B \times B$ on $(A, \mu|_A)$ and $(B, \nu|_B)$, are isomorphic. This notion is sometimes called *stable* orbit equivalence; however *weak* seems to be a more consistent terminology. Note that virtual isomorphism is a particular case of weakly orbit equivalence.

THEOREM A (strong orbit rigidity). *Let $\Gamma$ be a lattice in a simple, connected, noncompact Lie group $G$ with a finite center and $\mathbb{R} - \mathrm{rk}(G) \geq 2$. Let $\Gamma$ act ergodically by measure preserving transformations on a probability space $(X, \mathcal{B}, \mu)$. Assume that $X$ does not have $\Gamma$-equivariant measurable quotients of the form $\mathrm{Ad}\, G/\Delta$ where $\Delta \subset \mathrm{Ad}\, G$ is a lattice.*

*Then the $\Gamma$-action on $(X, \mu)$ is strongly orbitally rigid in the following sense*: *if $\Lambda$ is an arbitrary countable group which has an essentially free, f.m.p. ergodic action $(Y, \nu, \Lambda)$ which is weakly orbit equivalent to $(X, \mu, \Gamma)$, then $\Gamma$ and $\Lambda$ are virtually isomorphic groups and $(X, \mu, \Gamma)$ and $(Y, \nu, \Lambda)$ are virtually isomorphic actions.*

COROLLARY B. *The following actions $(X, \mu, \Gamma)$ are strongly orbitally rigid in the above sense*:

(1) *The natural action of $\mathrm{SL}_n(\mathbb{Z})$ on the n-torus $\mathbb{R}^n/\mathbb{Z}^n$, for $n > 2$.*

(2) *Let $\Gamma$ be a higher rank lattice as in Theorem A, and let $\Gamma \to K$ be a homomorphism into a compact group $K$ with a dense image. Then the action of $\Gamma$ on $K$ by translations is strongly orbitally rigid.*

(3) *Let $G \subset H$ be two simple, connected, noncompact Lie groups with finite centers where $\mathbb{R} - \mathrm{rk}(G) \geq 2$ and $G$ is not locally isomorphic to $H$. Let $\Gamma \subset G$ and $\Delta \subset H$ be lattices. Then the $\Gamma$-action on $H/\Delta$ is strongly orbitally rigid.*

Now consider the case of an ergodic action of a higher rank lattice $\Gamma \subset G$ which admits $G/\Delta$ quotients. The $\Gamma$-action on $G/\Delta$ itself is such an example, and it is known to be weakly orbit equivalent to the action of $\Delta$ on $G/\Gamma$ (see Example 2.10). The following theorem asserts that basically this is the only source of actions which are weakly orbit equivalent to an action of a higher rank lattice.

THEOREM C (orbit rigidity). *Let $\Gamma \subset G$ be a higher rank lattice as in Theorem A, and let $(X, \mu, \Gamma)$ be a f.m.p. ergodic action. Then, given any $\Gamma$-equivariant measurable quotient map $\pi : X \to \mathrm{Ad}\, G / \Delta$, where $\Delta \subset \mathrm{Ad}\, G$ is a lattice, there exists a canonically defined essentially free, f.m.p. ergodic action $(X_\pi, \mu_\pi, \Delta)$ of the lattice $\Delta$, which is weakly orbit equivalent to to the original action $(X, \mu, \Gamma)$.*

*Let $(Y, \nu, \Lambda)$ be any essentially free, f.m.p. and ergodic action of a countable group $\Lambda$, which is weakly orbit equivalent but not virtually isomorphic to $(X, \mu, \Gamma)$. Then there exists a lattice $\Delta \subset \mathrm{Ad}\, G$, and a measurable $\Gamma$-equivariant quotient map $\pi : X \to \mathrm{Ad}\, G / \Delta$, such that $\Lambda$ and $\Delta$ are virtually isomorphic groups, and $(Y, \nu, \Lambda)$ and $(X_\pi, \mu_\pi, \Delta)$ are virtually isomorphic actions.*

*Remark* 1.2. Any essentially free ergodic (measure preserving) action $(X, \mu, \Gamma)$ on a probability space $(X, \mu)$ gives rise to a factor von Neumann algebra (of type $\mathrm{II}_1$) $\mathbf{VN}(X, \Gamma)$ the so-called Murray-von Neumann group-space construction, which is a cross-product of the group algebra $\mathbf{VN}(\Gamma)$ with the abelian algebra $L^\infty(X, \mu)$ (the latter forms the so-called Cartan subalgebra of $\mathbf{VN}(X, \Gamma)$). The algebra $\mathbf{VN}(X, \Gamma)$ depends only on the orbit relation $(X, \mu, \mathcal{R}_\Gamma)$, and furthermore the relation $(X, \mu, \mathcal{R}_\Gamma)$ can be uniquely reconstructed (see [FM, Th. 1]) from the pair $(\mathbf{VN}(X, \Gamma), L^\infty(X, \mu))$. The above results show that for actions of higher rank lattices $\Gamma$, the pair consisting of the group space algebra $\mathbf{VN}(X, \Gamma)$ with its Cartan subalgebra $L^\infty(X, \mu)$ essentially determines the group $\Gamma$ and the action $(X, \mu, \Gamma)$ as in Theorems A and C above. This can be considered as supporting evidence to the far reaching Connes' conjecture that the von Neumann algebra $\mathbf{VN}(\Gamma)$ of a group $\Gamma$ with property T essentially determines the group $\Gamma$. Another evidence for this conjecture is the result of Cowling and Haagerup [CH] which in particular implies that the von Neumann algebras $\mathbf{VN}(\Gamma_i)$, $i = 1, 2$ of lattices $\Gamma_i \subset \mathrm{Sp}(n_i, 1)$ with $n_1 \neq n_2$ are not isomorphic.

Thus far we have considered equivalence relations $\mathcal{R} = \mathcal{R}_\Gamma$ on a probability space $(X, \mathcal{B}, \mu)$, which are generated by an *essentially free* action of some countable group $\Gamma$. One can also consider the more general setup of a countable relation of type $II_1$ (where $II_1$ corresponds to f.m.p. and nonatomic). For a detailed study of these and more general nonsingular equivalence relations we refer to the work of Feldman and Moore [FM]. In particular, at the very beginning of their analysis Feldman and Moore showed that any countable nonsingular relation on a probability space can be generated by a nonsingular action of some countable group. However, the following basic question remained open till now:

*Question* II ([FM, 77]). Is it true, that any countable nonsingular equivalence relation on a probability space can be generated by an *essentially free* action of some group?

Actually, the question formulated in [FM] addressed the possibility to generate freely a relation on a Borel space with no presence of a measure. This was answered in the negative by S. Adams in [Ad1]. However, the measure theoretic question above is more natural in the context of equivalence relations. The following theorem gives a negative solution to this problem:

THEOREM D. *The following are examples of equivalence relations of type $II_1$, which cannot be generated by an essentially free action of any group*:

(1) *Let $(X, \mu, \Gamma)$ be a strongly orbitally rigid action as in Theorem A and Corollary B, for example the natural action of $\mathrm{SL}_n(\mathbb{Z})$ on the torus $\mathbb{R}^n/\mathbb{Z}^n$. Let $A \subset X$ be a measurable subset with $\mu(A)$ being irrational. Then the restricted relation $(A, \mu_A, \mathcal{R}_\Gamma|_A)$ cannot be generated by an essentially free action of any group.*

(2) *Let $\Gamma$ be a higher rank lattice as in Theorem A. There exists a sequence of positive numbers $\{a_n\}$, such that for any f.m.p. ergodic $\Gamma$-action on a probability space $(X, \mathcal{B}, \mu)$ and any measurable subset $A \subset X$ such that $\mu(A)$ is not a rational multiple of some $a_n$, the restricted relation $(A, \mu_A, \mathcal{R}_\Gamma|_A)$ cannot be generated by an essentially free action of any group. The sequence $\{a_n\}$ consists of the covolume ratios $\mathrm{Haar}(G/\Gamma) : \mathrm{Haar}(G/\Delta)$, where $\Delta$ runs over all lattices in $G$.*

(3) *Let $G$ be as in Theorem A and let $G \subset H$ be a proper embedding of $G$ in a simple connected Lie group $H$. Consider $G$ acting on $H/\Delta$ where $\Delta$ is a lattice in $H$, and let $(Y, \mathcal{C}, \nu)$ be some measurable cross-section of this action with a $II_1$-relation $\mathcal{R}_G|_Y$ being defined by the intersection with the $G$-orbits (see Example 2.11). Then $(Y, \nu, \mathcal{R}_G|_Y)$ cannot be generated by an essentially free action of any group.*

The general idea behind these constructions is to modify an equivalence relation generated by a freely acting group with orbit equivalence rigidity properties in such a way that the original action does not generate the new relation any more, but that the rigidity properties of the action persist in some sense. Assuming that the new relation is generated by an (essentially) free action of some group, one uses the rigidity properties to show that the group and the action are the original ones, but that is impossible due to the modification made.

*Remark* 1.3. Construction (3) in Theorem D above is due to Zimmer, who had proved in [Zi3], that such $(Y, \nu, \mathcal{R}_G|_Y)$ cannot be essentially freely generated by a countable group, which admits a linear representation with an infinite image. It was then shown by Adams ([Ad2, 9.5]) that such $(Y, \nu, \mathcal{R}_G|_Y)$ cannot be essentially freely generated by an action of a word hyperbolic group. These results, however, left open the general question of an arbitrary group, which is resolved by showing that such freely acting group has to be virtually linear.

*Remark* 1.4. In the examples of type (1) and (2) in Theorem D the obstruction to being generated by an essentially free action lies in the "size" of the space: extending the space $(A, \mu_A)$ and the relation $\mathcal{R}_\Gamma|_A$ to the original $(X, \mu, \mathcal{R}_\Gamma)$ one achieves an orbit relation of an essentially free action. However, this "correction by extension" cannot be performed in the examples of type (3), since any "extension" $(Y_1, \nu_1, \mathcal{R}_1)$ of a given cross section $(Y, \nu, \mathcal{R}_G|_Y)$ can be realized as another cross section of the same $G$-action (Proposition 2.12).

Assertions (1) and (2) of Theorem D will follow from the proofs of Theorems A and C, while assertion (3) is a corollary of Theorem E below. Theorem E describes the situations where an ergodic $G$-action admits a measurable cross section $(Y, \nu, \mathcal{R}_G|_Y)$ of type $\mathrm{II}_1$, where $\mathcal{R}_G|_Y$ can be generated by an essentially free action of some group. Recall, that given any measure preserving, ergodic action $(X, \mu, \Gamma)$ of a lattice $\Gamma$ in a (general locally compact) group $G$, one can construct the corresponding induced $G$-action on $G/\Gamma \times X$ which is measure preserving and ergodic (cf. [Zi4, 4.2.21]). Such induced $G$-action has a fiber $(X, \Gamma)$ which might serve as a natural cross section of the $G$-action on $G/\Gamma \times X$, where the relation $\mathcal{R}_G|_X$ describing the intersections with $G$-orbits coincides with the original relation $\mathcal{R}_\Gamma$ on $X$ produced by the $\Gamma$-action. Hence starting with an essentially free action $(X, \mu, \Gamma)$ one obtains a $G$-action which admits cross sections generated by an essentially free action. Theorem E below states that in the case of higher rank simple Lie group these are the only examples.

THEOREM E. *Let $G$ be as in Theorem* A, *and let $G$ act measurably and ergodically on a probability space $(Z, \mathcal{D}, \eta)$ preserving the measure $\eta$. Assume that there exists a measurable cross section $(Y, \mathcal{C}, \nu)$ with $\mathrm{II}_1$-relation $\mathcal{R}_G|_Y$, which can be generated by an essentially free action of some countable group $\Lambda$. Then $(Y, \nu, \Lambda)$ is virtually isomorphic to some essentially free f.m.p. action $(Y_1, \nu_1, \Lambda_1)$ of a lattice $\Lambda_1 \subset \mathrm{Ad}\, G$, and the $G$-action on $(Z, \eta)$ is virtually isomorphic to the induced $\mathrm{Ad}\, G$-action on $\mathrm{Ad}\, G/\Lambda_1 \times Y_1$.*

Basically the same fact was obtained by Zimmer in [Zi3], under the assumption that $\Lambda$ admits a linear representation with an infinite image. The general case follows, therefore, as soon as one shows that whenever the relation $\mathcal{R}_G|_Y$ on a cross-section $(Y, \nu)$ is freely generated by some $\Lambda$, the latter is virtually linear. We remark also, that Theorem D.(2) shows that even for an induced $G$-action on $Z = G/\Gamma \times X$, for "most" cross sections $Y$ the relation $\mathcal{R}_G|_Y$ cannot be essentially freely generated by a group action.

Another natural question related to the induced actions is the following: let $(X, \mu, \Gamma)$ and $(Y, \nu, \Lambda)$ be two f.m.p. ergodic actions of two lattices $\Gamma, \Lambda \subset G$. Under what conditions are the induced $G$-action on $G/\Gamma \times X$ and the induced $G$-action on $G/\Lambda \times Y$ (virtually) isomorphic?

THEOREM F. *Let $G$ be as in Theorem* A *and let $(X, \mu, \Gamma)$ and $(Y, \nu, \Lambda)$ be measure preserving, ergodic actions of two lattices $\Gamma$ and $\Lambda$ in $G$. The induced $G$-actions on $G/\Gamma \times X$ and $G/\Lambda \times Y$ are virtually isomorphic if and only if the actions $(X, \mu, \Gamma)$ and $(Y, \nu, \Lambda)$ are weakly orbit equivalent. More precisely, either the actions $(X, \mu, \Gamma)$ and $(Y, \nu, \Lambda)$ themselves are virtually isomorphic, or there exist equivariant quotient maps $\pi_X : (X, \Gamma) \to (\mathrm{Ad}\, G/\mathrm{Ad}\, \Lambda, \Gamma)$ and $\pi_Y : (Y, \Lambda) \to (\mathrm{Ad}\, G/\mathrm{Ad}\, \Gamma, \Lambda)$ such that the weak orbit equivalence between $(X, \Gamma)$ and $(Y, \nu)$ factors through the standard weak orbit equivalence between $(\mathrm{Ad}\, G/\mathrm{Ad}\, \Lambda, \Gamma)$ and $(\mathrm{Ad}\, G/\mathrm{Ad}\, \Gamma, \Lambda)$ as in* (3.3).

*Organization of the paper.* Weak orbit equivalence, and more generally weak isomorphisms of equivalence relations, are discussed in Section 2. In Section 3 we consider the case of equivalence relations generated by essentially free actions of countable groups and describe the construction of the associated Gromov's measure equivalence couplings. With these preliminaries we prove all the above results in Section 4. The main tool in the proofs being the Measure Equivalence Rigidity Theorem from [Fu].

*Acknowledgments.* I would like to thank Robert Zimmer for his interest in this work and for many illuminating discussions. I would also like to thank Anatole Katok for his support and encouragement during the year that I spent at the Pennsylvania State University, as a Post Doctoral Fellow at the Center for Dynamical Systems.

## 2. Weak orbit equivalence

In this section we discuss a notion of weak isomorphism between countable $\mathrm{II}_1$-relations (or equivalently, weak orbit equivalence between f.m.p. group actions), and define a *compression constant* of a given weak isomorphism. Isomorphisms between $\mathrm{II}_1$-relations are characterized among all weak isomorphisms by the fact that they have compression constant one.

*Preliminaries.* First, let us recall the definition and some basic properties of $\mathrm{II}_1$-relations. A general countable nonsingular relation is an equivalence relation $\mathcal{R} \subset X \times X$ on a probability space $(X, \mathcal{B}, \mu)$, which forms a measurable subset $\mathcal{R} \in \mathcal{B} \times \mathcal{B}$, has at most countable equivalence classes $[x]_\mathcal{R}$ for $\mu$-a.e. $x \in X$, and satisfies the following property: for any zero measure set $A \subset X$ the saturation $\mathcal{R}(A) := \cup_{x \in A}[x]_\mathcal{R}$ of $A$, has also measure zero. We shall always assume that our equivalence relation $\mathcal{R}$ is ergodic in the sense that any measurable set which is a union of $\mathcal{R}$-classes has $\mu$-measure 0 or 1. An equivalence relation $\mathcal{R}$ on $(X, \mathcal{B}, \mu)$ is said to be of type $\mathrm{II}_1$, if there exists an $\mathcal{R}$-*invariant nonatomic probability* measure $\mu_0$ on $(X, \mathcal{B})$ which is equivalent to $\mu$: $\mu_0 \sim \mu$. Here being $\mathcal{R}$-invariant means that any Borel isomorphism $\phi : A \to B$ between measurable subsets $A, B \subseteq X$ with

$$(2.1) \qquad \mathrm{Graph}(\phi) := \{(x, \phi(x)) \in X \times X \mid x \in A\} \subset \mathcal{R}$$

satisfies $\mu_0(A) = \mu_0(B)$. In the sequel we consider only $\mathrm{II}_1$-relations $(X, \mu, \mathcal{R})$ and will always assume that $\mu$ itself is (the unique in its class) $\mathcal{R}$-invariant probability measure. It is easy to see that any ergodic, f.m.p. (not necessarily essentially free) action of a countable group $\Gamma$ on a nonatomic probability space $(X, \mu)$ generates an equivalence relation $\mathcal{R}_\Gamma$ of type $\mathrm{II}_1$. By [FM] any $\mathrm{II}_1$-relation can be presented in this form.

The collection of all partial Borel isomorphisms between subsets of $X$ whose graph is a subset of $\mathcal{R}$ as in (2.1) is denoted by $[[\mathcal{R}]]$. The subcollection of Borel isomorphisms $\phi : X \to X$ with $\mathrm{Graph}(\phi) \subset \mathcal{R}$ forms the so-called *full group* of $\mathcal{R}$, it is denoted by $[\mathcal{R}]$. If $\mathcal{R}$ is of type $\mathrm{II}_1$, as we assume, the elements of $[[\mathcal{R}]]$ and of $[\mathcal{R}]$ preserve the $\mathcal{R}$-invariant measure $\mu$. It is a standard (and not difficult) fact that the full group $[\mathcal{R}]$ acts transitively (mod 0) on the collection of measurable subsets of $X$ of a given measure.

Two nonsingular relations $\mathcal{R}$ on $(X, \mathcal{B}, \mu)$ and $\mathcal{S}$ on $(Y, \mathcal{C}, \nu)$ are said to be *isomorphic* (notation: $(X, \mu, \mathcal{R}) \cong (Y, \nu, \mathcal{S})$), if possibly after disregarding zero measure sets, there exists a Borel isomorphism $\theta : X \to Y$ with $\theta_*\mu \sim \nu$ and $\theta(\mathcal{R}) = \mathcal{S}$.

*Remark* 2.1. If $(X, \mu, \mathcal{R}) \cong (Y, \nu, \mathcal{S})$ are relations of type $\mathrm{II}_1$, then any isomorphism $\theta : X \to Y$ is measure preserving: $\theta_*\mu = \nu$. Indeed, any partial Borel isomorphism $\psi \in [[\mathcal{S}]]$ gives rise to a partial Borel isomorphism $\phi = \theta^{-1} \circ \psi \circ \theta \in [[\mathcal{R}]]$. Since $\phi_*\mu = \mu$ and $\psi_*\nu = \nu$, we have for $\nu$-a.e. $y$:

$$1 = \frac{d\psi_*\nu}{d\nu}(\psi(y)) = \frac{d\theta_*\mu}{d\nu}(\psi(y)) \cdot \frac{d\theta_*^{-1}\mu}{d\nu}(\theta^{-1}(y)) = \frac{d\theta_*\mu}{d\nu}(\psi(y)) : \frac{d\theta_*\mu}{d\nu}(y)$$

and therefore, $d\theta_*\mu/d\nu$ is a $\psi$-invariant function. Ergodicity of $\mathcal{S}$ implies that it is a.e. equal to a constant, which has to be $\nu(Y) : \mu(X) = 1$.

*Weak isomorphisms.* Given an equivalence relation $(X, \mu, \mathcal{R})$ and a positive measure subset $A \subseteq X$, one can consider the restricted equivalence relation $(A, \mu_A, \mathcal{R}|_A)$ defined by

$$\mu_A(E) := \mu(A \cap E)/\mu(A) \quad \text{and} \quad \mathcal{R}|_A := \mathcal{R} \cap (A \times A) .$$

*Definition* 2.2. Let $(X, \mu, R)$ and $(Y, \nu, S)$ be countable equivalence $II_1$-relations. We consider three equivalent definitions of *weak isomorphism* between the relations $(X, \mu, \mathcal{R})$ and $(Y, \nu, \mathcal{S})$, notation: $(X, \mu, \mathcal{R}) \simeq (Y, \nu, \mathcal{S})$:

1. There exist $A \subseteq X$, $\mu(A) > 0$ and $B \subseteq Y$, $\nu(B) > 0$, such that $(A, \mu_A, \mathcal{R}|_A) \cong (B, \nu_B, \mathcal{S}|_B)$.

2. There exists a $II_1$-relation $(Z, \eta, \mathcal{Q})$ and positive measure subsets $E$, $F \subseteq Z$ such that $(E, \eta_E, \mathcal{Q}|_E) \cong (X, \mu, \mathcal{R})$ and $(F, \eta_F, \mathcal{Q}|_F) \cong (Y, \nu, \mathcal{S})$.

3. There exist measurable maps $p : X \to Y$ and $q : Y \to X$, such that up to null sets

    (3a) $p_*\mu \prec \nu$ and $q_*\nu \prec \mu$,

    (3b) $p(\mathcal{R}) \subseteq \mathcal{S}$ and $q(\mathcal{S}) \subseteq \mathcal{R}$,

    (3c) $\text{Graph}(q \circ p) \subset \mathcal{R}$ and $\text{Graph}(p \circ q) \subset \mathcal{S}$.

We shall say that f.m.p. ergodic group actions $(X, \mu, \Gamma)$ and $(Y, \nu, \Lambda)$ are *weakly orbit equivalent* if the relations $(X, \mu, \mathcal{R}_\Gamma)$ and $(Y, \nu, \mathcal{R}_\Lambda)$ are weakly isomorphic.

In the introduction, we have used the first version of the definition, but in the sequel we shall mostly use the third one.

PROPOSITION 2.3. *Definitions* 2.2.1–2.2.3 *are equivalent.*

*Proof.* $1 \Rightarrow 2$. Let $\theta : A \to B$ be an isomorphism map between subsets $A \subseteq X$, $B \subseteq Y$. By Remark 2.1, $\theta$ is measure preserving: $\theta_*\mu_A = \nu_B$. Consider the union $Z = X \cup Y$, where $A$ is identified with $B$ via $\theta$, equipped with the finite positive measure $\eta_1$ defined by $\mu_1 = \mu(A)^{-1} \cdot \mu$ on $X$ and $\nu_1 = \nu(B)^{-1} \cdot \nu$ on $Y$, which agree on the identified sets. Normalizing $\eta_1$ to be a probability measure $\eta$ on $Z$, and considering the union $\mathcal{Q}$ of the equivalence relations $\mathcal{R}$ and $\mathcal{S}$, we obtain $(Z, \eta, \mathcal{Q})$ with the required properties.

$2 \Rightarrow 3$. Let $\theta_X : X \to E$ and $\theta_Y : Y \to F$ denote the isomorphisms. Since $E, F \subseteq Z$ have positive measure, and the relation $\mathcal{Q}$ is ergodic, one can construct measurable maps $\pi_E : Z \to E$ and $\pi_F : Z \to F$ with $\mathrm{Graph}(\pi_E)$, $\mathrm{Graph}(\pi_F) \subset \mathcal{Q}$ such that $\pi_E(x) = x$, $\pi_F(y) = y$ for $x \in E$, $y \in F$ (for example, let $\{\gamma_n\}_0^\infty$ be some enumeration of the elements of a group $\Gamma$ which generates $\mathcal{Q}$ with $\gamma_0 = e$; define $\pi_E(x) = \gamma_n \cdot x$ where $n$ is the first index with $\gamma_n \cdot x \in A$). Setting $p = \theta_Y^{-1} \circ \pi_F \circ \theta_X : X \to Y$ and $q = \theta_X^{-1} \circ \pi_E \circ \theta_Y : Y \to X$, one easily verifies the conditions (3a)–(3c). As in Remark 2.1, one can see that the constructed $(p, q)$, and generally any maps satisfying (3a)–(3c), are piecewise measure preserving in the following sense: there is a constant $c$, such that for any $E \subseteq X$ on which $p : E \to p(E)$ is one-to-one, we have $\nu(p(E)) = c \cdot \mu(E)$. Similar property holds for $q$ with $1/c$.

$3 \Rightarrow 1$. Let $p : X \to Y$ and $q : Y \to X$ satisfy (3a)–(3c). Let $B_1 := p(X) \subseteq Y$, $A_1 := q(B_1) \subseteq X$. Property (3a) implies $\nu(B_1) > 0$ and $\mu(A_1) > 0$. Let $B \subseteq B_1$ be some Borel cross section of $q : B_1 \to A_1$, and $A \subseteq X$ be some Borel cross section of $p : X \to B$. Then $\phi := q \circ p|_A$ is a measurable one-to-one map $\phi : A \to A_1$, and (3c) implies that $\phi \in [[\mathcal{R}]]$. Thus $\mu(A_1) = \mu(A) > 0$, so that $\theta := p|_A = q^{-1} \circ \phi|_A$ is a measure preserving isomorphism between positive measure subsets $A \subseteq X$ and $B \subseteq Y$, while property (3b) gives $\theta(\mathcal{R}|_A) = \mathcal{S}|_B$. □

*Definition* 2.4. Consider a weak orbit equivalence, given by maps $(p, q)$ as in Definition 2.2.3, between two *essentially free* actions $(X, \mu, \Gamma)$ and $(Y, \nu, \Lambda)$. Using property (3b) and the freeness assumption, we can define measurable functions $\alpha = \alpha_p : \Gamma \times X \to \Lambda$ and $\beta = \beta_q : \Lambda \times Y \to \Gamma$ by
$$p(\gamma \cdot x) = \alpha(\gamma, x) \cdot p(x), \qquad q(\lambda \cdot y) = \beta(\lambda, y) \cdot q(y).$$
The maps $\alpha$ and $\beta$ are measurable cocycles, i.e.
$$\alpha(\gamma_2 \gamma_1, x) = \alpha(\gamma_2, \gamma_1 \cdot x) \cdot \alpha(\gamma_1, x)$$
for $\mu$-a.e. $x$, and similarly for $\beta$. We shall say that $\alpha$ and $\beta$ are *weak orbit equivalent cocycles* (or wOE-cocycles) *associated* to $(p, q)$. If the action $(Y, \Lambda)$ is written from the right (as it will be in the sequel) we set

(2.2) $\qquad p(\gamma \cdot x) = p(x) \cdot \alpha(\gamma, x) \qquad q(y \cdot \lambda) = \beta(y, \lambda) \cdot q(y);$

here $\beta : Y \times \Lambda \to \Gamma$ becomes a right cocycle.

*Definition* 2.5. Let $(X, \mu, \mathcal{R})$ and $(Y, \nu, \mathcal{S})$ be weakly isomorphic $\mathrm{II}_1$-relations. Given a weak isomorphism $(X, \mu, \mathcal{R}) \simeq (Y, \nu, \mathcal{S})$ in the sense of Definitions 2.2.1–2.2.3, define a *compression constant* $c(X \simeq Y) \in (0, \infty)$ of this weak isomorphism as follows:

1. Given $\theta : (A, \mu_A, \mathcal{R}|_A) \cong (B, \nu_B, \mathcal{S}|_B)$ where $A \subseteq X, B \subseteq Y$, define:
$$c(X \simeq Y) = c(\theta) := \nu(B)/\mu(A).$$

2. Given $\theta_X : (E, \eta_E, \mathcal{Q}|_E) \cong (X, \mu, \mathcal{R})$ and $\theta_Y : (F, \eta_F, \mathcal{Q}|_F) \cong (Y, \nu, \mathcal{S})$ where $E, F \subseteq Z$, define:

$$c(X \simeq Y) = c(E, F) := \eta(E)/\eta(F).$$

3. Given $p : X \to Y$ and $q : Y \to X$, define

$$c(X \simeq Y) = c(p, q) := \nu(p(E))/\mu(E) = \nu(F)/\mu(q(F))$$

where $E \subseteq X$, $F \subseteq Y$ are any positive measure subsets, s.t. $p : E \to p(E)$ and $q : F \to q(F)$ are one-to-one.

Observe that the compression constant $c(X \simeq Y)$ depends not only on the relations $(X, \mu, \mathcal{R})$ and $(Y, \nu, \mathcal{S})$, but also on the weak isomorphism itself, given by $\theta$, $E \subset Z \supset F$ or $(p, q)$. However, all three types of weak isomorphisms, Definitions 2.2.1–2.2.3, which are related as in the proof of Proposition 2.3, have the same compression constant.

*Remark* 2.6. One can consider the set of all compression constants of weak isomorphisms of $(X, \mu, \mathcal{R})$ to itself. This set forms a subgroup, denote it $\mathrm{Comp}(\mathcal{R})$, of the multiplicative positive reals $\mathbb{R}_+^*$. This group is a subgroup of the so-called *fundamental group* of the von Neumann algebra associated with the relation $(X, \mu, \mathcal{R})$. For the (unique) amenable relation $\mathcal{R}_{\mathrm{amen}}$ one can show that $\mathrm{Comp}(\mathcal{R}_{\mathrm{amen}}) = \mathbb{R}_+^*$. However, for orbit relations produced by an action of a higher rank lattice $\Gamma$ it can be shown that $\mathrm{Comp}(\mathcal{R}_\Gamma) = \{1\}$. This latter property is shared by many nonamenable groups. We shall return to this subject elsewhere.

PROPOSITION 2.7. *Two* $\mathrm{II}_1$-*relations* $(X, \mu, \mathcal{R})$ *and* $(Y, \nu, \mathcal{S})$ *are isomorphic if and only if there exists a weak isomorphism with compression one*: $c(X \simeq Y) = 1$.

*Proof.* Obviously, any isomorphism of relations has compression constant one. Assume now that $(X, \mu, \mathcal{R})$ and $(Y, \nu, \mathcal{S})$ admit a weak isomorphism with $c(X \simeq Y) = 1$. Using definition (1) of 2.2 and 2.5, there exist subsets $A \subseteq X$, $B \subseteq Y$ with $\mu(A) = \nu(B) > 0$ and Borel isomorphism $\theta : A \to B$ with $\theta_* \mu_A = \nu_B$ and $\theta(\mathcal{R}|_A) = \mathcal{S}|_B$. The spaces $X$ and $Y$ can be decomposed into finite disjoint unions $X = \bigcup_0^k A_i$ and $Y = \bigcup_0^k B_i$ with $A_0 = A$, $B_0 = B$ and $\mu(A_i) = \nu(B_i) \leq \mu(A)$ for $1 \leq i \leq k$. Since the full group $[\mathcal{R}]$ acts transitively on subsets of equal measure, there exist $\phi_i \in [\mathcal{R}]$ and $\psi_i \in [\mathcal{S}]$, s.t. for $\phi_i(A_i) \subseteq A$ and $\psi_i(B_i) = \theta(A_i) \subseteq B$, for $0 \leq i \leq k$. Define a map $\Theta : X \to Y$ by setting $\Theta|_{A_i} := \psi_i^{-1} \circ \theta \circ \phi_i$. Then $\Theta$ is an isomorphism between $(X, \mu, \mathcal{R})$ and $(Y, \nu, \mathcal{S})$, which extends $\theta$; i.e. $\Theta|_A = \theta$. □

*Remark* 2.8. Similar arguments show that if $\theta_i : A_i \to B_i$, $i = 1, 2$, are two weak isomorphisms (in the sense of 2.2.1) of $(X, \mu, \mathcal{R})$ with $(Y, \nu, \mathcal{S})$, so that $(\theta_1 x, \theta_2 x) \in \mathcal{S}$ for a.e. $x \in X$, then $c(\theta_1) = c(\theta_2)$ if and only if there exist $\phi \in [\mathcal{R}]$, $\psi \in [\mathcal{S}]$ such that the maps $\theta_1$ and $\psi \circ \theta_2 \circ \phi$ coincide or one of the maps extends the other. With some abuse of terminology, we shall say that in this case the weak isomorphisms $\theta_1$ and $\theta_2$ are *conjugate* (with respect to the full groups). Similarly, two weak isomorphisms given by $(p_1, q_1)$ and $(p_2, q_2)$ as in 2.2.3 are *conjugate* (in the sense that $p_1|_A = \psi \circ p_2 \circ \phi|_A$ and $q_1|_B = \phi^{-1} \circ q_2 \circ \psi^{-1}|_B$ for some $\phi \in [\mathcal{R}]$, $\psi \in [\mathcal{S}]$ and positive measure subsets $A \subseteq X$, $B \subseteq Y$) if and only if $(p_1(x), p_2(x)) \in \mathcal{S}$ a.e. and $c(p_1, q_1) = c(p_2, q_2)$.

*Examples.* We shall now discuss some natural examples of weakly orbit equivalent actions.

*Example* 2.9. Let $\Gamma_1, \Gamma_2$ be virtually isomorphic groups and $(X_i, \mu_i, \Gamma_i)$, $i = 1, 2$, be virtually isomorphic actions as in Definition 1.1. Then there exists a natural weak orbit equivalence between the actions $(X_i, \mu_i, \Gamma_i)$ with the compression constant

$$c(X_1 \simeq X_2) = \frac{|N_2|}{|N_1|} \cdot \frac{[\Gamma_2' : \Gamma_2'']}{[\Gamma_1' : \Gamma_1'']} \cdot \frac{[X_2' : X_2'']}{[X_1' : X_1'']}$$

where $N_i \triangleleft \Gamma_i$ is finite, $\Gamma_i''$ is a finite index subgroup in $\Gamma_i' = \Gamma_i/N_i$, $(X_i', \mu_i')$ $= (X_i, \mu_i)/N_i$-quotient probability space, and $X_i'' \subseteq X'$ is one of at most $[\Gamma_i' : \Gamma_i'']$-many mutually isomorphic $\Gamma_i''$-ergodic components. We used the notation $[X_i' : X_i''] := \mu_i'(X_i')/\mu_i'(X_i'')$. Note that virtually isomorphic actions may have compression one, and hence being *orbit equivalent*, even when neither the actions nor the groups are isomorphic.

*Example* 2.10. Let $\Gamma$ and $\Lambda$ be lattices in a locally compact second countable (lcsc) group $G$. Then the $\Gamma$-action on $G/\Lambda$ and $\Lambda$-action on $G/\Gamma$ are weakly orbit equivalent with the compression constant

$$c(G/\Gamma \simeq G/\Lambda) = \mathrm{Haar}(G/\Gamma) : \mathrm{Haar}(G/\Lambda) .$$

Note that $G$ is unimodular (since it admits lattices) so that $g \mapsto g^{-1}$ gives an isomorphism of the left $\Lambda$-action on $G/\Gamma$ with the right $\Lambda$-action on $\Gamma\backslash G$. Choose right $\Lambda$- and left $\Gamma$- fundamental domains $X, Y \subset G$ which will represent the spaces $G/\Lambda$ and $\Gamma\backslash G$. Then the relations $(X, \mathcal{R}_\Gamma)$ and $(Y, \mathcal{R}_\Lambda)$ are given by the intersection of $X$ and $Y$ with the sets of the form $\Gamma g \Lambda$. All three definitions of weak isomorphism are easily verified now: assuming that $X \cap Y$ has positive Haar measure, one can take $A = B := X \cap Y$ and $\theta(x) = x$ in 2.2.1; or take $Z := X \cup Y$ for 2.2.2. Maps $(p, q)$ as in 2.2.3 can be defined by

$$p(x) := \Gamma x \cap Y \quad \text{and} \quad q(y) := X \cap y\Lambda .$$

The compression constant in all three options follows from the normalization to 1 of the Haar measure on the fundamental domains $X$, $Y$.

Another situation where weakly isomorphic relations appear naturally are measurable cross sections of group actions.

*Example* 2.11. Let $G$ be a locally compact second countable group, acting measurably on a probability space $(Z, \mathcal{D}, \eta)$. Assume that the action is measure preserving, properly ergodic (i.e. every orbit has zero measure) and locally free (i.e. the stabilizer of a.e. point is a discrete subgroups of $G$). Any such action admits a measurable cross section $Y \in \mathcal{D}$ which is equipped with a Borel structure $\mathcal{C}$ and a positive measure $\nu$, such that the the map $G \times Y \to Z$ given by $(g, y) \mapsto g \cdot y$ is measurable and is one-to-one measure preserving on a set $\bigcup_{y \in Y} U_y \times \{y\}$, where $U_y \subset G$ is a measurable family of open neighborhoods of the identity. Such cross section $(Y, \mathcal{C}, \nu)$ has a natural equivalence relation $\mathcal{R}_G|_Y$ defined by the intersections with the $G$-orbits:

$$\mathcal{R}_G|_Y := \{(y_1, y_2) \in Y \times Y \mid G \cdot y_1 = G \cdot y_2\}.$$

The measure $\nu$ is invariant and ergodic with respect to the countable relation $\mathcal{R}_G|_Y$, so that $(Y, \nu, \mathcal{R}_G|_Y)$ is either of type II$_1$ or II$_\infty$. In the latter case, passing to a subset $Y_0 \subset Y$ with $0 < \nu(Y_0) < \infty$, one obtains a new cross section $(Y_0, \mathcal{C}_0, \nu_0)$ with a II$_1$-relation $(Y_0, \nu_0, \mathcal{R}_G|_{Y_0})$.

The following easy proposition relates cross sections with the notion of weak isomorphism. The proof is analogous to the Kakutani Equivalence arguments. An excellent reference to the latter is the introduction of [ORW].

PROPOSITION 2.12. *Let $G$ be a lcsc group with a measurable, measure preserving, ergodic, locally free action on a probability space $(Z, \mathcal{D}, \eta)$. Then all measurable cross sections of this action $(Y, \mathcal{C}, \nu)$ with II$_1$-relations $\mathcal{R}_G|_Y$ are mutually weakly isomorphic.*

*If $G$ is not countable, then any II$_1$-relation $(Y', \nu', \mathcal{S})$, which is weakly isomorphic to one of the cross sections $(Y, \nu, \mathcal{R}_G|_Y)$, can be realized as a cross section $(Y', \nu', \mathcal{R}_G|_{Y'})$ of the same $G$-action on $(Z, \eta)$.*

## 3. Gromov's measure equivalence of groups

Let $(X, \mu, \Gamma)$ and $(Y, \nu, \Lambda)$ be *essentially free* actions which are orbit equivalent. By definition this means that there exists an isomorphism $\theta : X \to Y$ of orbit relations $\mathcal{R}_\Gamma$ and $\mathcal{R}_\Lambda$. In this setting the group actions are somewhat implicit, although a.e. orbit can be identified with the essentially freely acting group. Another aspect of orbit equivalence is the fact that the isomorphism $\theta$ is not canonically defined: for example each conjugate $\psi \circ \theta \circ \phi$ of $\theta$ (where $\phi \in [\mathcal{R}_\Gamma]$, $\psi \in [\mathcal{R}_\Lambda]$) gives an orbit equivalence $(X, \Gamma) \to (Y, \Lambda)$. The same applies to weak orbit equivalences.

However, in the case of (weakly) orbit equivalent actions which are essentially free one can construct a canonical object — the measure equivalent coupling, defined below — which has explicit group actions of $\Gamma$ and $\Lambda$ and which describes the entire conjugacy class of (weak) orbit equivalence maps between the given actions. The definition of measure equivalence (see below) was introduced by Gromov as a measure-theoretical analog of a "coarse geometric" notion of quasi-isometries between groups (see [Gr, 0.5.E] and [Fu]). The connection between measure equivalent couplings and orbit equivalent actions was presumably known to Gromov and was also observed by Zimmer, but apparently was not explicitly stated nor analyzed in detail in the literature. Therefore we do it here.

*Definition* 3.1 ([Gr, 0.5.E]). Two countable groups $\Gamma$ and $\Lambda$ are said to be *measure equivalent* if there exist commuting, measure preserving, essentially free actions of $\Gamma$ and $\Lambda$ on some infinite Lebesgue measure space $(\Omega, m)$, such that the action of each of the groups $\Gamma$ and $\Lambda$ admits a finite measure fundamental domain. The space $(\Omega, m)$ with the actions of $\Gamma$ and $\Lambda$ will be called a *measure equivalence coupling* of $\Gamma$ with $\Lambda$.

The basic example of measure equivalent groups are lattices in the same lcsc group $\Gamma, \Lambda \subset G$, where $(G, \text{Haar})$ with $\Gamma$ acting from the left, and $\Lambda$ from the right, gives a natural measure equivalence coupling of $\Gamma$ with $\Lambda$. Having this example in mind we shall use left- and right-action notations for actions in measure equivalence couplings.

Given a measure equivalent coupling $(\Omega, m)$ of some groups $\Gamma$ and $\Lambda$, one can choose fundamental domains $X, Y \subset \Omega$ for $\Omega/\Lambda$ and $\Gamma \backslash \Omega$. Given $\gamma \in \Gamma$ and a.e. $x \in X$, one can define a unique element $\alpha(\gamma, x) \in \Lambda$ by

(3.1) $$\gamma x \in X \, \alpha(\gamma, x) \ .$$

Define $\Gamma$-action on $X$ by: $\gamma \cdot x := \gamma x \alpha(\gamma, x)^{-1} \in X$. This action preserves the restriction $m|_X$ of $m$ to $X$. Identifying $X$ with the orbit space $\Omega/\Lambda$, we observe that the above f.m.p. $\Gamma$-action coincides with the natural $\Gamma$-action on $\Omega/\Lambda$. Similarly, define $\beta : Y \times \Lambda \to \Gamma$ and the right $\Lambda$-action on $Y$ by

$$y\lambda \in \beta(y, \lambda) \, Y, \qquad y \cdot \lambda := \beta(y, \lambda)^{-1} y \lambda \ .$$

This action coincides with the natural right $\Lambda$-action on $\Gamma \backslash \Omega$. Note, that $\alpha : \Gamma \times X \to \Lambda$ and $\beta : Y \times \Lambda \to \Gamma$ form left and right measurable cocycles, respectively:

$$\alpha(\gamma_2 \gamma_1, x) = \alpha(\gamma_2, \gamma_1 \cdot x) \, \alpha(\gamma_1, x) \qquad \beta(y, \lambda_1 \lambda_2) = \beta(y, \lambda_1) \, \beta(y \cdot \lambda_1, \lambda_2) \ .$$

We shall say that the cocycles $\alpha$ and $\beta$ are *measure equivalence cocycles* (or ME-cocycles) *associated* to the fundamental domains $X$ and $Y$ in the measure equivalent coupling $\Omega$. Different choices of fundamental domains $X$, $Y$ change the ME-cocycles $\alpha$, $\beta$, but do not change their measurable cohomology classes.

In the sequel, we shall be interested in *ergodic* measure equivalent couplings $\Omega$, which by definition means that $(\Omega/\Lambda, \Gamma)$ and $(\Gamma\backslash\Omega, \Lambda)$ are ergodic (cf. [Fu, 2.2]).

LEMMA 3.2. *Let $(\Omega, m)$ be some ergodic measure equivalent coupling of some groups $\Gamma$ and $\Lambda$. Let $X, Y \subset \Omega$ be some fundamental domains for $\Omega/\Lambda$ and $\Gamma\backslash\Omega$. Define maps $p_{X,Y} : X \to Y$ and $q_{X,Y} : Y \to X$ by*

(3.2) $$p_{X,Y}(x) := \Gamma x \cap Y, \qquad q_{X,Y}(y) := X \cap y\Lambda \ .$$

*Then $(p_{X,Y}, q_{X,Y})$ give a weak orbit equivalence (in the sense of 2.2.3) of the f.m.p. ergodic actions $(X, \Gamma)$ and $(Y, \Lambda)$ with the compression constant*

$$c(p_{X,Y}, q_{X,Y}) = m(Y) : m(X) \ .$$

*Moreover, the ME-cocycles $\alpha_{X,Y}$ and $\beta_{X,Y}$ associated to $X, Y \subset \Omega$ satisfy (2.2) and, therefore, they coincide with the wOE-cocycles $\alpha_{p_{X,Y}}$ and $\beta_{q_{X,Y}}$ associated to $(p_{X,Y}, q_{X,Y})$, provided that the actions $(X, \Gamma)$ and $(Y, \Lambda)$ are essentially free.*

*Proof.* The proof follows directly from the definitions. □

The following theorem enables to realize any weak orbit equivalence of essentially free f.m.p. actions in the above form.

THEOREM 3.3. *Let $(X, \mu, \Gamma)$ and $(Y, \nu, \Lambda)$ be ergodic f.m.p. actions of some countable groups $\Gamma$ and $\Lambda$. If the actions are weakly orbit equivalent and essentially free, then $\Gamma$ and $\Lambda$ are measure equivalent. Given a weak orbit equivalence $(p, q)$ as in 2.2.3, there exists an associated ergodic coupling $(\Omega_{p,q}, m_{p,q})$, such that*

$$(X, \Gamma) \cong (\Omega_{p,q}/\Lambda, \Gamma) \qquad (Y, \Lambda) \cong (\Gamma\backslash\Omega_{p,q}, \Lambda) \ .$$

*Moreover, $X$ and $Y$ can be realized as fundamental domains $\bar{X}, \bar{Y} \subset \Omega_{p,q}$ so that*

$$m_{p,q}(\bar{Y}) : m_{p,q}(\bar{X}) = c(p, q) \quad \text{and} \quad p = p_{\bar{X}, \bar{Y}} \quad q = q_{\bar{X}, \bar{Y}} \ .$$

*For any other weak orbit equivalence $(p', q')$ of $(X, \mu, \mathcal{R}_\Gamma)$ and $(Y, \nu, \mathcal{R}_\Lambda)$ satisfying $(p(x), p'(x)) \in \mathcal{R}_\Lambda$ for $\mu$-a.e. $x \in X$, the following are equivalent:*

- $(\Omega_{p',q'}, m_{p',q'}) \cong (\Omega_{p,q}, m_{p,q})$;

- $(p, q)$ and $(p', q')$ are conjugate (see 2.8);

- $c(p', q') = c(p, q)$.

*The case of orbit equivalence (rather than just a weak one) corresponds to a coupling $(\Omega, m)$ with fundamental domains of equal measure, in which case one can choose a common fundamental domain for both $\Gamma$- and $\Lambda$-actions.*

*Proof.* Let us write the $\Gamma$-action on $(X, \mu)$ from the left and the $\Lambda$-action on $(Y, \nu)$ from the right. The ME-cocycles

$$\alpha = \alpha_p : \Gamma \times X \to \Lambda \text{ and } \beta = \beta_q : Y \times \Lambda \to \Gamma$$

associated to $(p, q)$ as in (2.2) become respectively left and right cocycles. Hence we can consider the left $\Gamma$-action on the space $(\Omega_l, m_l) := (X, \mu) \times (\Lambda, m_\Lambda)$ and the right $\Lambda$-action on the space $(\Omega_r, m_r) := (\Gamma, m_\Gamma) \times (Y, \nu)$, defined by

$$\gamma(x, \lambda) = (\gamma \cdot x, \alpha(\gamma, x) \lambda), \qquad (\gamma, y) \lambda = (\gamma \beta(y, \lambda), y \cdot \lambda).$$

These are (infinite) measure preserving actions, which commute with the measure preserving right $\Lambda$- and left $\Gamma$-actions

$$(x, \lambda) \lambda_1 = (x, \lambda \lambda_1), \qquad \gamma_1 (\gamma, y) = (\gamma_1 \gamma, y).$$

The latter actions are essentially free and have finite measure fundamental domains $X \times \{e_\Lambda\}$ and $\{e_\Gamma\} \times Y$. We shall show that there exists a $\Gamma \times \Lambda$-equivariant isomorphism between $(\Omega_l, m_l)$ and $(\Omega_r, m_r)$.

Property (3c) of $p$ and $q$ states that there exist (unique) measurable assignments $x \to \gamma_x \in \Gamma$ and $y \to \lambda_y \in \Lambda$ such that $q(p(x)) = \gamma_x \cdot x$ and $p(q(y)) = y \cdot \lambda_y$. Consider the measurable maps $\Theta : (\Omega_l, m_l) \to (\Omega_r, m_r)$ and $\Theta' : (\Omega_r, m_r) \to (\Omega_l, m_l)$, defined by

$$\begin{aligned}\Theta(x, \lambda) &:= (\gamma_x \beta(p(x), \lambda), p(x) \cdot \lambda) \\ \Theta'(\gamma, y) &:= (\gamma \cdot q(y), \alpha(\gamma, q(y)) \lambda_y).\end{aligned}$$

Direct computation shows that $\Theta' = \Theta^{-1}$ and that $\Theta$ and $\Theta'$ are $\Gamma \times \Lambda$-equivariant. Since $p$ is (piecewise) measure preserving we have $\Theta_* m_l = c \cdot m_r$ where $c = c(X \simeq Y)$. Hence we obtain a measure equivalent coupling $(\Omega_{p,q}, m_{p,q}) := (\Omega_l, m_l) \cong (\Omega_r, c \cdot m_r)$ of $\Gamma$ with $\Lambda$ with fundamental domains

$$\bar{X} := X \times \{e_\Lambda\} \subset \Omega_l \cong \Omega_{p,q}, \qquad \bar{Y} := \{e_\Gamma\} \times Y \subset \Omega_r \cong \Omega_{p,q}.$$

We can naturally identify $(X, \mu)$ with $(\bar{X}, m_{p,q}|_{\bar{X}})$ and $(Y, c \cdot \nu)$ with $(\bar{Y}, m_{p,q}|_{\bar{Y}})$ as measure spaces with the actions of $\Gamma$ and $\Lambda$, respectively. It is easily seen that with this identification $(p, q)$ coincide with $(p_{\bar{X}, \bar{Y}}, q_{\bar{X}, \bar{Y}})$.

Observe, that any weak orbit equivalence $X \simeq Y$ given by $(p', q')$ with $c(p', q') = c(p, q)$ and $(p(x), p'(x)) \in \mathcal{R}_\Lambda$, can be realized as $(p_{\bar{X}', \bar{Y}'}, q_{\bar{X}', \bar{Y}'})$ for an appropriate choice of the fundamental domains $\bar{Y}'$ and $\bar{X}'$. So that the coupling $\Omega_{p', q'}$ will be isomorphic to $\Omega_{p,q}$ with the only difference being in the identification of $X$ and $Y$ with specific fundamental domains.

If $(X, \mu, \Gamma)$ and $(Y, \nu, \Lambda)$ are orbit equivalent, and $\theta = p = q^{-1} : X \to Y$ is the orbit equivalence relation isomorphism, then in the construction above, $\Theta(X \times \{e_\Lambda\}) = \{e_\Gamma\} \times Y$, so that $\Omega$ has one measurable subset $\bar{X} = \bar{Y}$, which is both $\Gamma$-and $\Lambda$- fundamental domain.

In the other direction, suppose that $(\Omega, m) = (\Omega_{p,q}, m_{p,q})$ is associated to a weak orbit equivalence $(p, q)$ of some essentially free $(X, \mu, \Gamma)$ and $(Y, \nu, \Lambda)$, so that $m(\Omega/\Lambda) = m(\Gamma \backslash \Omega)$. We can identify $X$ and $Y$ with some fundamental domains $X_0, Y_0 \subset \Omega$ with $m(X_0) = m(Y_0)$. We assumed that the $\Gamma$-action on $X \cong \Omega/\Lambda$ is ergodic. Hence $\Gamma \times \Lambda$-action on $(\Omega, m)$ is ergodic (of type $\mathrm{II}_\infty$). In this case, the action of the *full group* $[\Gamma \times \Lambda]$ (consisting of all Borel isomorphisms preserving the $\Gamma \times \Lambda$-orbits) is known to be transitive on sets of equal measure. This implies that there exist measurable (countable) partitions

$$X_0 = \bigcup_{\gamma \in \Gamma, \lambda \in \Lambda} X_{\gamma, \lambda} \quad \text{and} \quad Y_0 = \bigcup_{\gamma \in \Gamma, \lambda \in \Lambda} Y_{\gamma, \lambda} \quad \text{with} \quad Y_{\gamma, \lambda} = \gamma \, X_{\gamma, \lambda} \, \lambda.$$

Let $X_1 := \bigcup_{\gamma, \lambda} X_{\gamma, \lambda} \lambda$ and $Y_1 := \bigcup_{\gamma, \lambda} \gamma^{-1} Y_{\gamma, \lambda}$. Note that, being formed by piecewise $\Lambda$-translations of $X_0$, the set $X_1$ gives a $\Lambda$-fundamental domain; similarly $Y_1$ gives a $\Gamma$-fundamental domain. But this is the same set: $X_1 = Y_1$. □

*Quotients of measure equivalent couplings.* The measure equivalence point of view on orbit equivalence will be useful in the proofs of our results. In particular, consider the following situation. Let $(\Omega_0, m_0)$ and $(\Omega_1, m_1)$ be two (ergodic) measure equivalent couplings of some groups $\Gamma$ and $\Lambda$, and assume that $(\Omega_0, m_0)$ is a $\Gamma \times \Lambda$-equivariant quotient of $(\Omega_1, m_1)$, i.e. there exists a measurable map $\Phi : \Omega_1 \to \Omega_0$ with

$$\Phi_* m_1 = m_0, \quad \text{and} \quad \Phi(\gamma \, \omega_1 \, \lambda) = \gamma \, \Phi(\omega_1) \, \lambda, \quad (\gamma \in \Gamma, \ \lambda \in \Lambda, \ \omega_1 \in \Omega_1).$$

Let $X_0$, $Y_0$ be some fundamental domains for $\Omega_0/\Lambda$ and $\Gamma \backslash \Omega_0$. Then the preimages

$$X_1 := \Phi^{-1}(X_0) \qquad Y_1 := \Phi^{-1}(Y_0)$$

form fundamental domains for $\Omega_1/\Lambda$ and $\Gamma \backslash \Omega_1$. Let us denote by $(p_i, q_i)$ the maps $(p_{X_i, Y_i}, q_{X_i, Y_i})$ associated to this choice of the fundamental domains. By Lemma 3.2 $(p_i, q_i)$ give weak orbit equivalences $(X_i, \Gamma) \simeq (Y_i, \Lambda)$ for $i = 0, 1$. Since the the disintegration of $m_1$ with respect to $m_0 = \Phi_* m_1$ consists of finite measures, which are therefore probability measures, the compression constants coincide

$$c(p_1, q_1) = c(p_0, q_0).$$

Moreover, the restrictions of $\Phi$ to $X_1$ and $Y_1$ give equivariant quotient maps $\pi_X : (X_1, \Gamma) \to (X_0, \Gamma)$ and $\pi_Y : (Y_1, \Lambda) \to (Y_0, \Lambda)$, which commute with $(p_i, q_i)$:

(3.3)
$$\begin{array}{ccc} (X_1, \Gamma) & \xrightarrow{p_1} & (Y_1, \Lambda) \\ \pi_X \downarrow & & \downarrow \pi_Y \\ (X_0, \Gamma) & \xrightarrow{p_0} & (Y_0, \Lambda) \end{array} \quad \text{and} \quad \begin{array}{ccc} (X_1, \Gamma) & \xleftarrow{q_1} & (Y_1, \Lambda) \\ \pi_X \downarrow & & \downarrow \pi_Y \\ (X_0, \Gamma) & \xleftarrow{q_0} & (Y_0, \Lambda) \end{array}.$$

Consider the ME-cocycles $\alpha_i : \Gamma \times X_i \to \Lambda$ and $\beta_i : Y_i \times \Lambda \to \Gamma$, $i = 0, 1$ which are defined by the maps $(p_i, q_i)$ as in (3.1). Commutativity of the diagram (3.3) implies that for a.e. $x_1 \in X_1$, $y_1 \in Y_1$ and all $\gamma \in \Gamma$, $\lambda \in \Lambda$:

(3.4) $\qquad \alpha_1(\gamma, x_1) = \alpha_0(\gamma, \pi_X(x_1)) \quad \beta_1(y_1, \lambda) = \beta_0(\pi_Y(y_1), \lambda)$.

If the actions $(X_i, \Gamma)$ and $(Y_i, \Lambda)$ are essentially free, ME-cocycles $\alpha_i$ and $\beta_i$ coincide with the wOE-cocycles (defined by (2.2)), which describe the correspondence of $\Gamma$-orbits to $\Lambda$-orbits under $(p_i, q_i)$. Thus, viewing $\alpha_i$ and $\beta_i$ as wOE-cocycles, we conclude from (3.4) that the weak orbit equivalence $(X_1, \Gamma) \simeq (Y_1, \Lambda)$ basically takes place on the level of the quotients $(X_0, \Gamma) \simeq (Y_0, \Lambda)$ and is just extended without further rearrangements in the fibers. We summarize this discussion in the following:

PROPOSITION 3.4 (reduction to quotients). *Consider a weak orbit equivalence $(p, q)$ between two essentially free f.m.p. actions $(X_1, \Gamma) \simeq (Y_1, \Lambda)$, which has the property that the associated measure equivalence coupling $\Omega$ has a $\Gamma \times \Lambda$-equivariant quotient coupling $\Omega_0$. Then $(p, q)$ can be replaced by a conjugate (see 2.8) weak orbit equivalence $(p_1, q_1)$ which factors through the weak orbit equivalence $(p_0, q_0)$ between the $\Gamma$-action on $X_0 := \Omega_0 / \Lambda$ and the $\Lambda$-action on $Y_0 := \Gamma \backslash \Omega_0$, so that (3.3) and (3.4) hold. The compression constants of the original weak orbit equivalence and the constructed ones are the same*:

$$c(p, q) = c(p_1, q_1) = c(p_0, q_0).$$

Now consider the question of reconstructing $(Y_1, \Lambda)$ from the quotient map $(X_1, \Gamma) \to (X_0, \Gamma)$ and the weak orbit equivalence $(Y_0, \Lambda) \simeq (X_0, \Gamma)$.

PROPOSITION 3.5 (lifting of weak orbit equivalence). *Let $(X_1, \mu_1, \Gamma)$, $(X_0, \mu_0, \Gamma)$ and $(Y_0, \nu_0, \Lambda)$ be essentially free f.m.p. ergodic actions of some countable groups $\Gamma$ and $\Lambda$. Assume that there exists a $\Gamma$-equivariant quotient map $\pi_X : X_1 \to X_0$ and a weak orbit equivalence $(p_0, q_0)$ between $(X_0, \mu_0, \Gamma)$ and $(Y_0, \nu_0, \Lambda)$. Then there exist: an essentially free f.m.p. ergodic action of $\Lambda$ on some probability space $(Y_1, \nu_1)$, a weak orbit equivalence $(p_1, q_1)$ of $(X_1, \Gamma)$ and $(Y_1, \Lambda)$ and a $\Lambda$-equivariant quotient map $\pi_Y : Y_1 \to Y_0$, such that the diagram (3.3) commutes and (3.4) holds. The action $(Y_1, \nu_1, \Lambda)$ and the maps $(p_1, q_1)$ and $\pi_Y$ are uniquely determined by (3.3)*.

We remark that if one replaces the condition that (3.3) commutes by a weaker condition that $(Y_1, \nu_1, \Lambda)$ is weakly orbit equivalent to $(X_1, \mu_1, \Gamma)$ and have the given $(Y_0, \nu_0, \Lambda)$ as an equivariant quotient, then the uniqueness part of the proposition does not hold, in general.

*Proof of Proposition* 3.5. Let $(\Omega_0, m_0)$ be the measure equivalent coupling of $\Gamma$ and $\Lambda$ associated with $(p_0, q_0)$, so that $X_0$ and $Y_0$ are realized as fundamental domains $\bar{X}_0, \bar{Y}_0$, so that that $p_0(x) = \Gamma x \cap Y_0$ and $q_0(y) = X_0 \cap y\Lambda$ (Theorem 3.3). We shall construct another $\Gamma, \Lambda$ coupling $(\Omega_1, m_1)$ for which $(X_1, \Gamma) \cong (\Omega_1/\Lambda, \Gamma)$, and will take $Y_1 := \Gamma \backslash \Omega_1$.

Let $\alpha_0 = \alpha_{p_0} : \Gamma \times X_0 \to \Lambda$ be the weak orbit equivalence cocycle associated with $p_0$ as in (2.2). Consider the space $\Omega_1 := X_1 \times \Lambda$ with the measure $m_1 := \mu_1 \times m_\Lambda$ where $m_\Lambda$ is the counting measure. $\Lambda$ acts from the right on $\Omega_1$ (just in the $\Lambda$ coordinate) and $\Gamma$ acts from the left by

$$\gamma : (x, \lambda) \mapsto (\gamma \cdot x, \alpha_0(\gamma, \pi_X(x))\lambda) .$$

Consider the map $\Phi : \Omega_1 \to \Omega_0$ defined by $\Phi(x, \lambda) = \pi_X(x)\lambda$, where $\pi_X(x) \in X_0$ is viewed as a point in $\Omega_0$. Then $\Phi$ is $\Gamma$-equivariant and $\Phi_* m_1 = m_0$. Setting

$$Y_1 := \Phi^{-1}(Y_0), \quad \pi_Y := \Phi|_{Y_1}, \quad p_1(x) := \Gamma x \cap Y_1, \quad q_1(y) := X_1 \cap y\Lambda ,$$

we are in the situation of Proposition 3.4 and the existence part of the claim follows.

As for the uniqueness, observe that the construction of $\Omega_1$ is completely determined by $(p_0, q_0)$ and $\pi_X$. This implies the uniqueness of the action $(Y_1, \Lambda)$. The realization of $Y_1$ as $\Phi^{-1}(Y_0)$ and the definition of the maps $p_1, q_1, \pi_Y$ are also dictated by the requirement that the diagram (3.3) should commute. □

## 4. Proofs of the main results

The main tool in the proofs is the following result from [Fu]:

THEOREM 4.1 (measure equivalence rigidity of higher rank lattices). *Let $\Gamma$ be a lattice in a simple, connected, noncompact Lie group $G$ with finite center and $\mathbb{R} - \mathrm{rk}(G) \geq 2$. Let $\Lambda$ be an arbitrary countable group, which is measure equivalent to $\Gamma$. Then $\Lambda$ is virtually isomorphic to a lattice in $\mathrm{Ad}\, G$, i.e. there exists a finite index subgroup $\Lambda' \subseteq \Lambda$ and a homomorphism $\rho : \Lambda' \to \mathrm{Ad}\, G$ with $N = \mathrm{Ker}(\rho)$ being finite and $\rho(\Lambda')$ being a lattice in $\mathrm{Ad}\, G$.*

*Moreover, if $(\Omega, m)$ is an ergodic coupling of $\Gamma$ with $\Lambda$, then there exists a unique measurable map $\Phi : \Omega \to \mathrm{Ad}\, G$ satisfying*

$$\Phi(\gamma^{-1} x \lambda) = \mathrm{Ad}\,(\gamma)^{-1}\, \Phi(x)\, \rho(\lambda)$$

*for $m$-a.e. $x \in \Omega$ and all $\gamma \in \Gamma$, $\lambda \in \Lambda'$. The $\Phi$-projection $\Phi_* m$ of the measure $m$ to $\operatorname{Ad} G$ is either the Haar measure on $\operatorname{Ad} G$, or an atomic measure. The fibers of the disintegration of $m$ with respect to $\Phi_* m$ are probability measures.*

The statement of the Main Theorem (3.1) in [Fu] uses $\operatorname{Aut}(\operatorname{Ad} G)$ and $\Lambda$ itself (rather than $\operatorname{Ad} G$ and $\Lambda'$), however the above formulation is easily deduced applying the finite-to-one map $\operatorname{Aut}(\operatorname{Ad} G) \to \operatorname{Ad} G$. For our applications we have also stated the result in the case of ergodic measure equivalence coupling $(\Omega, m)$, where the mixture of atomic and Haar measures cannot occur for $\Phi_* m$.

The above theorem should be interpreted as follows: up to virtual isomorphisms, lattices in a given higher rank $G$ form a single class of measure equivalent groups. Moreover any ergodic measure equivalent coupling of two groups $\Gamma$ and $\Lambda$ from this class has a (canonically defined) natural quotient coupling, which is either $\operatorname{Ad} G$ with its Haar measure or some atomic space. We shall use both parts of this statement in the following proofs.

*The Main Step.* Let $\Gamma$ be a lattice in $G$ as in Theorem A, acting ergodically by measure preserving transformations on a non atomic Lebesgue probability space $(X, \mu)$. Let $(Y, \nu, \Lambda)$ be an essentially free f.m.p. action of some countable group $\Lambda$, and assume that $(X, \mu, \Gamma)$ and $(Y, \nu, \Lambda)$ are weakly orbit equivalent. Passing, if necessary, to the action $(X_1, \mu_1, \Gamma_1)$ of $\Gamma_1 := \Gamma/\Gamma \cap Z(G)$, we obtain an essentially free action ([SZ]). Let us write the $\Lambda$-action on $Y$ from the right. Applying Theorem 3.3, we construct a measure equivalent coupling $(\Omega, m)$ of $\Gamma_1$ with $\Lambda$, such that the left $\Gamma_1$-action on $X_1$ and on $\Omega/\Lambda$, and the right $\Lambda$-actions on $Y$ and on $\Gamma\backslash\Omega$ are isomorphic.

Now, Theorem 4.1 implies that $\Lambda$ is virtually isomorphic to a lattice $\rho(\Lambda')$ in $\operatorname{Ad} G$. Replacing $\Lambda$ by $\Lambda_1 := \rho(\Lambda')$ and passing to the corresponding action $(Y_1, \nu_1, \Lambda_1)$ and to the measure equivalence coupling $(\Omega_1, m_1)$ we obtain the situation of two lattices $\Gamma_1 \subset \operatorname{Ad} G \supset \Lambda_1$ and a measure equivalent coupling $(\Omega_1, m_1)$ with a (unique) $\Gamma_1 \times \Lambda_1$-equivariant measurable map

(4.1) $\quad \Phi_1 : \Omega_1 \to G, \qquad \Phi_1(\gamma \omega \lambda) = \gamma \Phi_1(\omega) \lambda \qquad (\gamma \in \Gamma_1,\ \lambda \in \Lambda_1,\ \omega \in \Omega_1)$

where $m_0 := \Phi_{1*} m_1$ is the Haar measure on $\operatorname{Ad} G$ or is an atomic measure. The actions $(X_1, \Gamma_1) \cong (\Omega_1/\Lambda_1, \Gamma_1)$ and $(Y_1, \Lambda_1) \cong (\Gamma_1\backslash\Omega_1, \Lambda_1)$ are essentially free, weakly orbit equivalent to each other, and are virtually isomorphic to the original actions $(X, \Gamma)$ and $(Y, \Lambda)$.

Fix some $\Lambda_1$- and $\Gamma_1$- fundamental domains $X_0, Y_0 \subseteq \operatorname{Ad} G$, let $\mu_0$ and $\nu_0$ be the normalized restrictions $m_0|_{X_0}$ and $m_0|_{Y_0}$, and consider the f.m.p. actions $(X_0, \mu_0, \Gamma_1)$ and $(Y_0, \nu_0, \Lambda_1)$ with the weak orbit equivalence $(p_0, q_0)$ given by

$$p_0(x) := \Gamma_1 x \cap Y_0 \qquad q_0(y) := X_0 \cap y \Lambda_1\ .$$

By Proposition 3.4, there exist equivariant maps $\pi_X : X_1 \to X_0$ and $\pi_Y : Y_1 \to Y_0$ and a new weak orbit equivalence $(p_1, q_1)$, which is conjugate to the previous one, so that the diagram (3.3) commutes and (3.4) holds.

We shall deduce Theorems A and C from the two cases: (i) $m_0 = \Phi_{1*}m_1$ is atomic, or (ii) $m_0 = \Phi_{1*}m_1$ coincides with the Haar measure on $\operatorname{Ad} G$, respectively.

*Proof of Theorem* A. It is assumed in the theorem that $(X, \mu, \Gamma)$ does not have equivariant quotients of the form $(\operatorname{Ad} G/\Delta, \Gamma)$. The passage to $(X_1, \Gamma_1) = (X, \Gamma)/(\Gamma \cap Z(G))$ has not changed this property. Hence, in this case $m_0 = \Phi_{1*}m_1$ is an atomic measure on $\operatorname{Ad} G$, and therefore $\mu_0$ and $\nu_0$ on $X_0$ and $Y_0$ are *finite* atomic measures. This implies that $\Gamma_1$ and $\Lambda_1$ are commensurable; more precisely, for some (any) atom at $g$ in the support of $m_0$, the group

$$\Gamma_2 = \Gamma_1 \cap g\Lambda_1 g^{-1}$$

is of finite index in both $\Gamma_1$ and in $g\Lambda_1 g^{-1}$. We can assume that $X_0, Y_0 \subset \operatorname{Ad} G$ were chosen so that $g \in X_0 \cap Y_0$. Let $\Lambda_2 := g^{-1}\Gamma_2 g \subseteq \Lambda_1$ and consider $X_2 = Y_2 := \Phi_1^{-1}(\{g\})$ which is both $\Gamma_2$- and $\Lambda_2$-invariant subset in $\Omega_1$. Note that $X_1, Y_1$ as the subsets $\Phi^{-1}(X_0), \Phi^{-1}(Y_0) \subset \Omega_1$ consist of finitely many $\Gamma_1$- and $\Lambda_1$-translations of the common set $X_2 = Y_2$. The weak orbit equivalence maps $(p_1, q_1)$ were defined by

$$p_1(x) = \Gamma_1 x \cap Y_1, \qquad q_1(y) = X_1 \cap y\Lambda_1 .$$

Therefore, on the set $X_2 = Y_2 \subseteq X_1 \cap Y_1$ we have $p_1 = q_1 = \operatorname{Id}$. Thus the cocycles $\alpha_1 : \Gamma_1 \times X_1 \to \Lambda_1$ and $\beta_1 : Y_1 \times \Lambda_1 \to \Gamma_1$, associated to $(p_1, q_1)$ (both as ME-cocycles and as wOE-cocycles), have the property that the restrictions $\alpha_1|_{\Gamma_2 \times X_2}$ and $\beta_1|_{Y_2 \times \Lambda_2}$ do not depend on the space variable (we use (3.4) and the fact that $g$ is fixed by the $\Gamma_2$-action on $X_0$, and by the $\Lambda_2$-action on $Y_0$). This means that these restricted cocycles are given by homomorphisms $a : \Gamma_2 \to \Lambda_2$ and $b : \Lambda_2 \to \Gamma_2$ and thus $a = b^{-1}$ is an isomorphism of groups, and the actions $(X_2, \Gamma_2)$ and $(Y_2, \Lambda_2)$ are isomorphic with respect to this group isomorphism $a = b^{-1}$. We conclude that the pairs of actions $(X_1, \Gamma_1)$ and $(Y_1, \Lambda_1)$, and hence also the original $(X, \Gamma)$ and $(Y, \Lambda)$, are virtually isomorphic actions. □

*Proof of Corollary* B. Note that the systems $(X, \mu, \Gamma)$ listed in the corollary do not have equivariant quotients of the form $(\operatorname{Ad} G/\Delta, \operatorname{Haar}, \Gamma)$:

(1) Observe that for $\gamma \in \operatorname{SL}_n(\mathbb{R})$ which has at least one eigenvalue not on the unit circle, the entropy $h(\gamma, \mathbb{T}^n)$ of the action on the torus is strictly less than the entropy $h(\gamma, \operatorname{PSL}_n(\mathbb{R})/\Delta)$, since these entropies can be

computed by the formulas

$$h(\gamma, \mathbb{T}^n) = \sum_{\lambda_i \geq 0} \lambda_i, \qquad h(\gamma, \mathrm{PSL}_n(\mathbb{R})/\Delta) = \sum_{\lambda_i \geq \lambda_j} \lambda_i - \lambda_j$$

where $\exp \lambda_i$ is the absolute value of the $i^{\mathrm{th}}$ eigenvalue of $\gamma$.

(2) Obvious from an even simpler entropy argument.

(3) Follows from Witte [Wi] Corollary 5.6′.

*Proof of Theorem* C. We fix an ergodic f.m.p. action $(X, \mu, \Gamma)$ of a higher rank lattice $\Gamma \subset G$. Let $\pi$ be an equivariant projection $\pi : X \to \mathrm{Ad}\, G/\Delta$ where $\Delta$ is a lattice in $\mathrm{Ad}\, G$. Dividing, if necessary, the group $\Gamma$ and the action $(X, \mu, \Gamma)$ by the finite group $\Gamma \cap Z(G)$, we can reduce the discussion to the case where $\Gamma = \mathrm{Ad}\,\Gamma$ is a lattice in $\mathrm{Ad}\, G$. Choose some fundamental domains $X_0, Y_0 \subset \mathrm{Ad}\, G$ for $\mathrm{Ad}\, G/\Delta$ and $\Gamma \backslash \mathrm{Ad}\, G$, and let $(p_0, q_0)$ be the associated weak orbit equivalence maps between $(X_0, \Gamma)$ and $(Y_0, \Delta)$, which are essentially free f.m.p. ergodic actions. In Proposition 3.5 we have constructed an essentially free f.m.p. $\Delta$-action $(Y_1, \nu_1, \Delta)$, weak orbit equivalence maps $(p_1, q_1)$ and a $\Delta$-equivariant quotient map $\pi : Y_1 \to Y_0$, so that diagram (3.3) commutes. Such an action $(Y_1, \nu_1, \Delta)$ is uniquely determined, and it was denoted by $(X_\pi, \mu_\pi, \Delta)$ in the statement of the theorem.

Now let $(Y, \nu, \Lambda)$ be any essentially free f.m.p. action of some countable $\Lambda$, which is weakly orbit equivalent but not virtually isomorphic to $(X, \mu, \Gamma)$. Applying Theorem 4.1 and The Main Step above, we obtain a lattice $\Lambda_1 \subset \mathrm{Ad}\, G$ which is virtually isomorphic to $\Lambda$, and the action $(Y_1, \nu_1, \Lambda_1)$ which is virtually isomorphic to $(Y, \nu, \Lambda)$. Moreover, $(Y_1, \nu_1, \Lambda_1)$ is weakly orbit equivalent, by $(p_1, q_1)$, to $(X_1, \mu_1, \Gamma_1) := (X, \mu, \Gamma)/\Gamma \cap Z(G)$, and this equivalence factors through a weak orbit equivalence $(p_0, q_0)$ of the equivariant quotients

$$\begin{array}{rcl}(X_1, \mu_1) \xrightarrow{\pi_X} (X_0, \mu_0) & \cong & (\mathrm{Ad}\, G, m_0)/\Lambda_1 \\ (Y_1, \nu_1) \xrightarrow{\pi_Y} (Y_0, \nu_0) & \cong & \Gamma_1 \backslash (\mathrm{Ad}\, G, m_0)\ .\end{array}$$

Since $(X, \mu, \Gamma)$ and $(Y, \nu, \Lambda)$ are assumed not to be virtually isomorphic, $m_0$ has to be the Haar measure. In particular, we obtain a $\Gamma_1$-equivariant quotient map $\pi_X : (X_1, \mu_1) \to (\mathrm{Ad}\, G/\Lambda_1, \mathrm{Haar})$ which gives rise to the $\Gamma$-equivariant quotient map:

$$\pi : (X, \mu) \longrightarrow (X_1, \mu_1) \xrightarrow{\pi_X} (\mathrm{Ad}\, G/\Lambda_1, \mathrm{Haar})\ .$$

Let $(X_\pi, \mu_\pi, \Lambda_1)$ be the $\Lambda_1$-action constructed above. Then the actions $(Y_1, \nu_1, \Lambda_1)$ and $(X_\pi, \mu_\pi, \Lambda_1)$ both fit in the commuting diagram (3.3), and by the uniqueness part of Proposition 3.5, are isomorphic. This completes the proof of Theorem C. □

*Proof of Theorem* D. To prove that examples (1) and (2) of Theorem D give a negative solution to Feldman-Moore's problem (question II), we shall use orbit equivalence rigidity described in Theorems A and C together with the compression constants. Let $(X, \mu, \Gamma)$ be an essentially free f.m.p. ergodic action of a higher rank lattice $\Gamma \subset G$. And suppose that a restricted relation $\mathcal{R}_\Gamma|_A$ on some positive measure subset $A \subset X$ is isomorphic to a relation $\mathcal{R}_\Lambda$ on $(Y, \nu)$, where $(Y, \nu, \Lambda)$ is an *essentially free* f.m.p. action of some countable group $\Lambda$. Such isomorphism $\theta : A \to Y$ forms a *weak* orbit equivalence between essentially free f.m.p. actions $(X, \mu, \Gamma)$ and $(Y, \nu, \Lambda)$. Such $\theta$ can be replaced by some $(p, q)$ as in Proposition 2.3 with the same compression constant

$$c(p, q) = c(\theta) = 1/\mu(A) .$$

Following the proofs of Theorems A and C, we observe that applying some virtual isomorphisms (The Main Step), we obtained a weak orbit equivalence $(p'_1, q'_1)$ between $(X_1, \mu_1, \Gamma_1)$ and $(Y_1, \nu_1, \Lambda_1)$ and quotient actions $(X_0, \mu_0, \Gamma_1)$ and $(Y_0, \nu_0, \Lambda_1)$. The weak orbit equivalence $(p'_1, q'_1)$ was replaced then by a conjugate one $(p_1, q_1)$ which factors through $(p_0, q_0)$ in the sense of (3.3). The crucial observation is that the compression constants satisfy

$$r \cdot c(p, q) = c(p'_1, q'_1) = c(p_1, q_1) = c(p_0, q_0)$$

where the rational multiple $r \in \mathbb{Q}$ corresponds to the virtual isomorphisms used in the Main Step. For actions $(X, \mu, \Gamma)$ satisfying Theorem A, the constant $c(p_0, q_0)$ is rational, and hence so is $c(p, q) = 1/\mu(A)$. For more general actions $(X, \mu, \Gamma)$ of a higher rank lattice $\Gamma \subset G$, we have

$$c(p_0, q_0) = \mathrm{Haar}(\mathrm{Ad}\, G/\mathrm{Ad}\,\Gamma) : \mathrm{Haar}(\mathrm{Ad}\, G/\Lambda_1);$$

hence $\mu(A)$ has to be a rational multiple of a covolume ratio of lattices in $G$. It is well known that the covolumes of lattices form a countable set, and, in fact, there are only finitely many covolumes in $[0, M]$ for any $M < \infty$. This proves assertions (1) and (2) of Theorem D. Assertion (3) will be deduced from Theorem E by the observation that an induced $G$-action on $G/\Lambda_1 \times Y$ cannot be isomorphic to a $G$-action on $H/\Delta$, where $G \subset H$ is a proper subgroup in a simple Lie group $H$. For the latter see Witte [Wi] Corollary 5.6′.

*Proof of Theorem* E. We start with an ergodic f.m.p. action of a higher rank $G$ on some $(Z, \eta)$ which admits a cross section $(Y, \nu)$ with the $\mathrm{II}_1$-relation $\mathcal{R}_G|_Y$ being generated by an essentially free action of some countable $\Lambda$. We shall write this $\Lambda$-action on $(Y, \nu)$ from the right. Since the cross section $\mathcal{R}_G|_Y$ is assumed to be of type $\mathrm{II}_1$(and not atomic) the $G$-action has to be properly ergodic, i.e. every orbit has zero measure, and therefore by [SZ], it is essentially free. Hence we can define a measurable map $\alpha : Y \times \Lambda \to G$ by

$$y \cdot \lambda = \alpha(y, \lambda)\, y, \qquad (\lambda \in \Lambda,\ g \in G,\ y \in Y) .$$

Observe that $\alpha$ is a right cocycle. Choose some lattice $\Gamma$ in $G$. We shall construct a measure equivalence coupling $(\Omega, m)$ of $\Gamma$ with $\Lambda$ by the following: let $\Omega := G \times Y$ and $m := \text{Haar}_G \times \nu$, and consider the actions

$$\gamma(g,y) := (\gamma g, y) \qquad (g,y)\lambda := (g\alpha(y,\lambda), y \cdot \lambda) \ .$$

Actually we have a m.p. essentially free $G$-action on $\Omega$, where $\Gamma$ acts as a subgroup. $\Gamma$- and $\Lambda$- actions are commuting, measure preserving and *essentially free*, where the freeness of the $\Lambda$-actions follows from the freeness assumption in the $Y$-coordinate. The $\Gamma$-action has a finite measure fundamental domain $G/\Gamma \times Y$. We claim that the $\Lambda$-action also has a finite measure fundamental domain. Indeed, let $\pi : \Omega \to Z$ be the map defined by $\pi(g,y) := g \cdot y$. Since $Y$ is a cross section, $\pi$ is a countable-to-one measurable map which is piecewise measure preserving. More precisely, any Borel map $\sigma : Z \to \Omega$ satisfying $\pi \circ \sigma(z) = z$, is measure preserving. One easily verifies that the image $\sigma(Z) \subset \Omega$ of such $\sigma$ forms a fundamental domain for the $\Lambda$-action on $\Omega$.

We now apply Theorem 4.1 to conclude that $\Lambda$ is virtually isomorphic to a lattice $\Lambda_1 \subset \text{Ad}\, G$, i.e. $\Lambda$ has a finite index subgroup $\Lambda'$ so that $\Lambda_1 \cong \Lambda'/N$ with $N$ being a finite normal subgroup. Consider the left $G$-action on $(Z', \eta') := (\Omega, m)/\Lambda'$ which is a finite cover of the $G$-action on $(Z, \eta) = (\Omega, m)/\Lambda$. The action $(Z', \eta', G)$ has a cross section $(Y', \nu')$, consisting of $[\Lambda : \Lambda']$-translations of $(Y, \nu)$. Moreover, there is an action of $\Lambda'$ on $(Y', \nu')$ such that $\mathcal{R}_G|_{Y'} = \mathcal{R}_{\Lambda'}$. Dividing $(\Omega, m)$ by the actions of the finite center $Z(G)$ from the left, and by the $N$-action from the right, we obtain a left $\text{Ad}\, G$-action and the right $\Lambda_1$-action on

$$(\Omega_1, m_1) := Z(G)\backslash(\Omega,m)/N = \text{Ad}\, G \times (Z(G)\backslash Y'/N) \ .$$

The $\text{Ad}\, G$-action on $(Z_1, \eta_1) := (\Omega_1, m_1)/\Lambda_1$ is a finite-to-one quotient of $(Z', \eta', G)$, and it has a cross section $(Y_1, \nu_1) = Z(G)\backslash(Y', \nu')/N$ which is equipped with a $\Lambda_1$-action satisfying $\mathcal{R}_{\text{Ad}\, G}|_{Y_1} = \mathcal{R}_{\Lambda_1}$. Obviously, $(Y, \nu, \Lambda)$ and $(Y_1, \nu_1, \Lambda_1)$ are virtually isomorphic.

Since $\Lambda_1$ is a lattice in $\text{Ad}\, G$, one can construct the induced $\text{Ad}\, G$-action on

$$(\text{Ad}\, G/\Lambda_1 \times Y_1, \text{Haar} \times \nu_1) \ .$$

It has a cross section $(Y_1, \nu_1)$ which is intersected by $G$-orbits along $\mathcal{R}_{\Lambda_1}$. The same happens with the $\text{Ad}\, G$-action on $(Z_1, \eta_1)$. Hence the two $\text{Ad}\, G$-actions are orbit equivalent, and therefore, by (Zimmer [Zi2], or [Zi4, 5.2.1]) these $\text{Ad}\, G$-actions are isomorphic. □

*Proof of Theorem* F. Let $(X, \mu, \Gamma)$ and $(Y, \nu, \Lambda)$ be two f.m.p. ergodic actions of lattices $\Gamma, \Lambda$ in a higher rank $G$. Assume that the actions are weakly orbit equivalent, i.e. there exists a map $X \supset A \xrightarrow{\theta} B \subset Y$ which is an isomorphism of the relations $\theta : (A, \mu_A, \mathcal{R}_\Gamma|_A) \cong (B, \nu_B, \mathcal{R}_\Lambda|_B)$. Consider the induced $G$-actions on $G/\Gamma \times X$ and $G/\Lambda \times Y$. These actions have cross sections $(X, \mu)$ and $(Y, \nu)$ with the relations $\mathcal{R}_G|_X = \mathcal{R}_\Gamma$ and $\mathcal{R}_G|_Y = \mathcal{R}_\Lambda$, respectively. Since $(A, \mu_A, \mathcal{R}_\Gamma|_A)$ and $(B, \nu_B, \mathcal{R}_\Lambda|_B)$ form cross sections as well, we conclude that the $G$-actions on $G/\Gamma \times X$ and on $G/\Lambda \times Y$ are orbit equivalent. Using superrigidity for cocycles Zimmer proved ([Zi2], or [Zi4, 5.2.1]), that in such case the $G$-actions are virtually isomorphic.

On the other hand, if the induced $G$-actions on $G/\Gamma \times X$ and on $G/\Lambda \times Y$ are virtually isomorphic, division by the finite group $Z(G)$ gives an isomorphism of the $\mathrm{Ad}\, G$-actions

$$(\mathrm{Ad}\, G/\Gamma_1 \times X_1, \mathrm{Ad}\, G) \cong (\mathrm{Ad}\, G/\Lambda_1 \times Y_1, \mathrm{Ad}\, G)\ .$$

Since $(X_1, \mu_1, \mathcal{R}_{\Gamma_1})$ and $(Y_1, \nu_1, \mathcal{R}_{\Lambda_1})$ are cross sections of the same (up to isomorphism) $\mathrm{Ad}\, G$-action, they are weakly isomorphic relations (Prop. 2.12) and, therefore, $(X, \mu, \Gamma)$ and $(Y, \nu, \Lambda)$ are weakly orbit equivalent actions. We are in the situation, discussed at the beginning of this section and the assertions of the theorem follow from this discussion. Actually for the current case, where $\Lambda$ is already known to form a lattice in $G$, we do not need the whole strength of Theorem 4.1, but rather the first step of its proof, namely [Fu, 4.1]. □

*Current address*: University of Illinois at Chicago, Chicago IL
*E-mail address*: furman@math.uic.edu